\documentclass[11pt]{article}
\usepackage{amsmath,latexsym,epsf,mathabx}
\usepackage[all]{xy}
\usepackage{mathrsfs}
\usepackage{soul}
\usepackage{color}
\begin{document}

\newcommand{\nc}{\newcommand}
\newtheorem{Th}{Theorem}[section]
\newtheorem{Def}[Th]{Definition}
\newtheorem{Lem}[Th]{Lemma}
\newtheorem{Pro}[Th]{Proposition}
\newtheorem{Cor}[Th]{Corollary}
\newtheorem{Rem}[Th]{Remark}
\newtheorem{Exm}[Th]{Example}
\newtheorem{Sc}[Th]{}
\def\Pf#1{{\noindent\bf Proof}.\setcounter{equation}{0}}
\def\bskip#1{{ \vskip 20pt }\setcounter{equation}{0}}
\def\sskip#1{{ \vskip 5pt }\setcounter{equation}{0}}
\def\mskip#1{{ \vskip 10pt }\setcounter{equation}{0}}
\def\bg#1{\begin{#1}\setcounter{equation}{0}}
\def\ed#1{\end{#1}\setcounter{equation}{0}}

\soulregister\cite7 
\soulregister\citep7 
\soulregister\citet7 
\soulregister\ref7


\title{\bf  Balanced pairs and tilting modules in recollement
\thanks{Supported by the National Natural Science Foundation of China (Grants No.11801004) and the Startup Foundation for Introducing Talent of AHPU (Grant No.2017YQQ016). }
}
\smallskip
\author{ Peiyu Zhang and Dajun Liu\\ 
\footnotesize ~E-mail:~zhangpy@ahpu.edu.cn, ldjnnu2017004@163.com; \\
\footnotesize School of Mathematics and Physics,  Anhui Polytechnic University, Wuhu, China.  \\
\\ Jiaqun Wei\\ 
\footnotesize ~E-mail:~weijiaqun@njnu.edu.cn;\\
\footnotesize  Department of Mathematics, Northwest Normal University, Lanzhou, China;\\
\footnotesize  School of Mathematics Science, Nanjing Normal University, Nanjing, China.
}

\date{}
\maketitle
\baselineskip 15pt
%
%
\begin{abstract}
\vskip 10pt%
In this paper, firstly, we mainly study the relationship of balanced pairs among three
Abelian categories in a recollement. As an application of admissible balanced pairs, we introduce the notion
of the relative tilting modules, and give a characterization of relative tilting modules, which similar to Bazzoni characterization of $n$-tilting modules \cite{BS}.
Finally, we mainly consider the relationship of relative tilting modules in a recollement.

\mskip\




\end{abstract}
%
\vskip 30pt

\section{Introduction}

A recollement of Abelian categories is special algebra structure consisting of three Abelian categories and six functors.
It has applications in many aspects of algebra, for example, representation theory, ring theory, geometry etc.
Recollements of Abelian categories were used by Cline, Parshall and Scott to study module categories of finite dimensional algebras over a field in \cite{PS}.
Recently, C. Psaroudakis and J. Vitória established a correspondence between recollements of abelian categories up to equivalence and certain TTF-triples in \cite{PC1}.
In addition, the notion of Recollements of abelian categories and recollements of triangulated categories are similar. For instance, the author \cite{CJM} constructed a recollement
of abelian categories from a recollement of triangulated categories, generalizing a result of \cite{LYN}. In this paper, we only consider the recollement of Abelian categories.
There are many examples of recollements, see \cite{PC}.

The author \cite{CHEN} introduced the notion of balanced pairs in the study of triangle-equivalence between two homotopy categories of complexes.
Furthermore, Chen proved that there exists a triangle-equivalence between the homotopy category of Gorenstein projective modules and the homotopy
category of Gorenstein injective modules when $R$ is a left-Gorenstein ring. Zheng construct a recollement of triangulated categories through
a complete balanced pair in \cite{zyf}.
In section 2, we mainly consider whether the recollement of Abelian categories has the property of preserve balance pairs.
i.e., how to construct a balance pair from given balance pair(s) in a recollement, see Theorem \ref{result1} and Theorem \ref{result2}.

The tilting theory is well known, the research of tilting theory is abundant.
Angeleri-H\"{u}gel and Coelho \cite{AC}, Bazzoni \cite{BS} and Wei\cite{WEI} considered tilting modules of projective dimension $\leq n$, that is, $n$-tilting modules.
As an application of admissible balanced pair, we give the notion of the relative tilting modules in section 3 of this paper.
Bazzoni's characterization of the $n$-tilting modules is extended the relative tilting modules, see Theorem\ref{mainresult}.
Finally, we mainly study the relationship of the relative tilting modules among three
Abelian categories in a recollement, see Theorem \ref{TH1} and Theorem \ref{TH2}.
%
%
%
\hskip 18pt

\section{Balanced pairs and Recollement}
\vskip 10pt

Let the functor $F$: ${\mathscr A}\longrightarrow {\mathscr B}$ and the functor $G$: ${\mathscr B}\longrightarrow {\mathscr A}$, where both $\mathscr A$ and $\mathscr B$
are Abelian categories. we said $(F$, $G)$ to be an adjoint pair,
if there is a isomorphism $\sigma_{X,Y}$: $\mathrm{Hom}_{{\mathscr B}}(FX, Y)\cong\mathrm{Hom}_{{\mathscr A}}(X, GY)$ for any $X\in {\mathscr A}$ and $Y\in {\mathscr B}$.
As we all know, an adjoint pair $(F$, $G)$ induces two natural transformations $\eta$: $\mathrm{Id}_{{\mathscr A}}\longrightarrow GF$ and
$\varepsilon$: $FG\longrightarrow \mathrm{Id}_{{\mathscr B}}$. Set $\eta_{X}:=\sigma_{X,~F(X)}(\mathrm{Id}_{F(X)})$: $X\longrightarrow GF(X)$ and
$\varepsilon_{Y}:=\eta^{-1}_{G(Y),~Y}(\mathrm{Id}_{G(Y)})$: $FG(Y)\longrightarrow Y$ for any $X\in {\mathscr A}$, $Y\in {\mathscr B}$, which called
the unit and counit of the adjunction, respectively, such that $\sigma_{X,Y}(f)=G(f)\eta_{X}$ and $\sigma^{-1}_{X,Y}(g)=\varepsilon_{Y} F(g) $
for any $f\in \mathrm{Hom}_{\mathscr A}(F(X),Y)$, $g\in \mathrm{Hom}_{\mathscr B}(X,G(Y))$.

\mskip\
The following conclusion is well known.

\bg{Lem}\label{adjoindpro}
Let $(F$, $G)$ is an adjoint pair.

$(1)$ We have that $Id_{F(X)}=\varepsilon_{F(X)}F(\eta_{X})$ and $Id_{G(Y)}=G(\varepsilon_{Y})\eta_{G(Y)}$ for any $X\in \mathscr A$, $Y\in \mathscr B$;

$(2)$ $F~(resp.~G)$ is fully faithful if and only if $\eta~(resp.~\varepsilon)$ is an isomorphism.
\ed{Lem}

\bg{Def}$\mathrm{\cite{FP,PC1}}$\label{RE}%
A recollement of an abelian category ${\mathscr A}$ by abelian categories ${\mathscr A'}$ and
${\mathscr A''}$, denoted by $R(\mathscr A', ~\mathscr A, ~\mathscr A'')$, is a diagram of additive functors as follows, satisfying the conditions below.

$$\xymatrix{
{\mathscr A^\prime} \ar[rrr]|{\ i_{\ast}}&&& {\mathscr A}
\ar@/^/@<2ex>[lll]|{i^{!}}\ar[rrr]|{\
j^{\ast}}\ar@/_/@<-2ex>[lll]|{i^{\ast}} &&& {\mathscr A''}
\ar@/_/@<-2ex>[lll]|{j_{!}}\ar@/^/@<2ex>[lll]|{j_{\ast}} }$$

$(i)$ ($i^{\ast}$, $i_{\ast}$, $i^{!}$) and ($j_{!}$, $j^{\ast}$, $j_{\ast}$) are adjoint triples;

$(ii)$ The functors $i_{\ast}$, $j_{!}$, and $j_{\ast}$ are fully faithful;

$(iii)$ $\mathrm{Im}i_{\ast}$ = $\mathrm{Ker}j^{\ast}$.

\ed{Def}

\bg{Exm}$\mathrm{\cite{PC,PC1}}$
Let A be a ring, $e\in A$ satisfying $e^{2}=e$. There is a recollement of modules category $R(\mathrm{Mod}-A/AeA,~\mathrm{Mod}-A, ~\mathrm{Mod}-eAe)$ as follows,
which is said to be induced by the idempotent $e$.
$$\xymatrix{
\mathrm{Mod}-A/AeA \ar[rrr]^{inc}&&& \mathrm{Mod}-A
\ar@/^/@<2ex>[lll]^{\mathrm{Hom}_{A}(A/AeA,-)}\ar[rrr]^{\mathrm{Hom}_{A}(eA,-)}\ar@/_/@<-2ex>[lll]_{-\otimes_{A}A/AeA} &&& \mathrm{Mod}-eAe
\ar@/_/@<-2ex>[lll]_{-\otimes_{eAe}eA}\ar@/^/@<2ex>[lll]^{\mathrm{Hom}_{eAe}(Ae,-)}}$$
\ed{Exm}

In fact, there is many examples of recollements, see Example 2.8-2.13 in \cite{PC}.
Next, we collect some properties of recollements, which is very useful in the sequel \cite{FP,mh,PC,PSS,PC1}.

\bg{Pro}\label{propR}%

Let $R(\mathscr A', ~\mathscr A, ~\mathscr A'')$ be a recollement of abelian categories. Then we have the following properties.

$(1)$ $i^{\ast}j_{!}=0$ and $i^{!}j_{\ast}=0$;

$(2)$  $i^{\ast}$ and $j_{!}$ are right exact, $i^{!}$ and $j_{\ast}$ are left exact, $i_{\ast}$ and $j^{\ast}$ are exact;

$(3)$ These natural transformations $i^{\ast}i_{\ast}\longrightarrow \mathrm{Id}_{\mathscr A'}$, $\mathrm{Id}_{\mathscr A'}\longrightarrow i^{!}i_{\ast}$,
$j^{\ast}j_{\ast}\longrightarrow \mathrm{Id}_{\mathscr A''}$, $\mathrm{Id}_{\mathscr A''}\longrightarrow j^{\ast}j_{!}$ are natural isomorphisms.

$(4)$ If $i^{\ast}$ is exact, then $i^{!}j_{!}=0$ and $j_{!}$ is exact; If $i^{!}$ is exact, then $i^{\ast}j_{\ast}=0$ and $j_{\ast}$ is exact.

\ed{Pro}

A left $\mathscr C$-resolution of $M$ is a complex $\cdots \longrightarrow C_{2}\longrightarrow C_{1}\longrightarrow C_{0}\longrightarrow M$
with $C_{i}\in \mathscr C$ for $i\geq0$ such that it is acyclic by applying the functor $\mathrm{Hom}_{\mathscr C}(C,-)$ for each $C\in\mathscr C$.
We denote sometimes the left $\mathscr C$-resolution of $M$ by $C^{\bullet}\longrightarrow M$,
where $C^{\bullet}=:$ $\cdots \longrightarrow C_{2}\longrightarrow C_{1}\longrightarrow C_{0}\longrightarrow 0$
is the deleted left $\mathscr C$-resolution of $M$.
The left $\mathscr C$-dimension of $M$, written $\mathscr C$-dim $M$, is defined as inf $\{ n|$
there is a exact sequence $0\longrightarrow C_{n}\longrightarrow \cdots \longrightarrow C_{1}\longrightarrow C_{0}\longrightarrow M$
with each $C_{i} \in \mathcal{C}$ $\}$.
If no such an integer exists, then $\mathscr C$-dim $M=\infty$.
Define the global $\mathscr C$-dimension $\mathscr C$-dim $\mathscr A$ to be the supreme of the $\mathscr C$-resolution dimensions of all the objects
in $\mathscr A$.

Recall that a subcategory $\mathscr C$ of $\mathscr A$ is said to be contravariantly
finite, if for any $A\in \mathscr A$, it has a right $\mathscr C$-approximation \cite{AR}, i.e., there is
a homomorphism $f$ : $C\to A$ for some $C\in \mathscr C$ such
that $\mathrm{Hom}_{\mathscr A}(C',f)$ is surjective for any $C'\in \mathscr C$.
Dually, we have the definition of covariantly finite subcategory.

\bg{Def}$\mathrm{\cite{CHEN}}$\label{BP}%
A pair $({\mathscr X}$, ${\mathscr Y})$ of additive subcategories in ${\mathscr A}$ is called a balanced pair if the following
conditions are satisfied:

$(1)$ the subcategory ${\mathscr X}$ is contravariantly finite and ${\mathscr Y}$ is covariantly finite;

$(2)$ for each object M, there is an ${\mathscr X}$-resolution $X^{\bullet}\longrightarrow M$ such that it is acyclic by applying the
functors $\mathrm{Hom}_{{\mathscr A}}(-,~Y )$ for all $Y\in{\mathscr Y}$;

$(3)$ for each object N, there is a ${\mathscr Y}$-coresolution $N\longrightarrow Y^{\bullet}$ such that it is acyclic by applying the
functors $\mathrm{Hom}_{{\mathscr A}}(X,~-)$ for all $X\in{\mathscr X}$.
\ed{Def}

By the proposition 2.2 in \cite{CHEN}, if ${\mathscr X}$ (resp., ${\mathscr Y}$) be a contravariantly finite (resp., covariantly finite) subcategory,
then the pair $({\mathscr X}$, ${\mathscr Y})$ is balanced if and only if the class of right ${\mathscr X}$-acyclic complexes
coincides with the class of left ${\mathscr Y}$-acyclic complexes.

\bg{Exm}\label{eg1}
$(1)$ Let R be an ring. It is easily see that $(\mathrm{Proj} R,~\mathrm{Inj} R)$ is a balanced pair in $\mathrm{Mod}$-R, where
$\mathrm{Proj} R$ $(respectively, \mathrm{Inj} R)$ is all projective (respectively, injective) modules of $\mathrm{Mod}$-R.

$(2)$ According to \cite{gj,kh}, over a commutative noetherian ring with finite Krull dimension, $(\mathcal{GP},~\mathcal{GP}^{\bot})$ and $(^{\bot}\mathcal{GI},~\mathcal{GI})$ are complete and hereditary cotorsion pairs \cite{EJ}, where $\mathcal{GP}$ (respectively, $\mathcal{GI}$) is
the subcategory of all Gorenstein projective (respectively, Gorenstein injective) modules \cite{Holm}. A commutative noetherian ring $R$ called virtually
Gorenstein \cite{za}, if it has finite Krull dimension and $\mathcal{GP}^{\bot}=^{\bot}\mathcal{GI}$. Therefore, in this case, $(\mathcal{GP},~\mathcal{GP}^{\bot}=^{\bot}\mathcal{GI},~\mathcal{GI})$ is a complete and hereditary cotorsion triple. By Proposition 2.6 in \cite{CHEN}, we can obtain that $(\mathcal{GP},~\mathcal{GI})$ is a balanced pair.

$(3)$ $($See \cite[Example 8.3.2]{EJ}$)$ Let R be a ring and $\mathrm{Mod}$-R the category of right R-modules, and let
$\mathcal{PP}(R)$ and $\mathcal{PI}(R)$ be the subcategories of $\mathrm{Mod}$-R consisting of pure projective modules and pure
injective modules respectively. Then $(\mathcal{PP}(R),~\mathcal{PI}(R))$ is a balanced pair in $\mathrm{Mod}$-R.
\ed{Exm}

\mskip\
A contravariantly finite subcategory $\mathscr C$ of $\mathscr A$ is admissible provided that each right
$\mathscr C$-approximation is surjective. It is equivalent to that any right $\mathscr C$-acyclic complex is acyclic. Indeed,
if $\mathscr C$ is admissible, we consider the right $\mathscr C$-acyclic complex
$\xymatrix{0\ar[r]&M\ar[r]^{f} & N\ar[r]^{g}& L\ar[r] &0}$
with $f$ injective and $g$ surjective. Suppose that $C_{M}^{\bullet}\longrightarrow M$ and $C_{L}^{\bullet}\longrightarrow L$ are left $\mathscr C$-resolution
of $M$ and $L$, respectively. We can easily obtain the following commutative diagram:
$$\xymatrix{
0\ar[r]&M\ar[r] & N\ar[r]& L\ar[r]&0\\
0\ar[r]&C_{M}^{\bullet}\ar[r]\ar[u]^{\varepsilon_{M}} &C_{N}^{\bullet}\ar[r]\ar[u]^{\varepsilon_{N}} &C_{L}^{\bullet}\ar[r]\ar[u]^{\varepsilon_{L}}&0}$$
where $C_{N}^{\bullet}=C_{M}^{\bullet}\bigoplus C_{L}^{\bullet}$. Since $\mathscr C$ is admissible, $\varepsilon_{M}$ and
$\varepsilon_{L}$ are quasi-isomorphic, thus, $\varepsilon_{N}$ is also quasi-isomorphic. Consequently,
the first row is acyclic since the second row is acyclic. On the contrary, for any $M\in \mathscr A$, there is an exact sequence
$\xymatrix{0\ar[r]&K\ar[r]&C\ar[r]^{\alpha}& M}$ with $\alpha$ right $\mathcal{C}$-approximation. Since the complex
$$\xymatrix{0\ar[r]&\mathrm{Hom}_{\mathscr A}(C^{'},K)\ar[r]&\mathrm{Hom}_{\mathscr A}(C^{'},C)\ar[r]& \mathrm{Hom}_{\mathscr A}(C^{'},M)\ar[r] &0}$$
is exact for any $C^{'}\in \mathscr C$, $\alpha$ is surjective.

Recalled that a balanced pair $({\mathscr X}$, ${\mathscr Y})$ is admissible in the Abelian category ${\mathscr A}$, if ${\mathscr X}$ is admissible.
By the corollary 2.3 in \cite{CHEN}, a balanced pair $({\mathscr X}$, ${\mathscr Y})$ is admissible if and only if ${\mathscr Y}$ is admissible.
If the $\mathscr X$-dim $\mathscr A$ is finite, then the balanced pair $({\mathscr X}$, ${\mathscr Y})$ is of finite dimensions.

\mskip\

The following lemma is the version of Lemma 3.1 in \cite{LYN} in Abelian category.

\bg{Lem}\label{YL1}
Let F and G be two functors, $\xymatrix{{\mathscr A} \ar@<1ex>[r]^F& {\mathscr B} \ar@<1ex>[l]^{G} }$. If ${\mathscr C}\subseteq {\mathscr A}$ is a
contravariantly finite subcategory of ${\mathscr A}$, then $F({\mathscr C})$ is a contravariantly finite subcategory of ${\mathscr B}$.
Moreover, if G is faithful and ${\mathscr C}$ is admissible, then $F({\mathscr C})$ is also admissible.

\ed{Lem}

\Pf. For any $B\in\mathscr B$, there is a right $\mathscr C$-approximation $f$: $C\longrightarrow G(B)$ of $G(B)$ since ${\mathscr C}$ is a contravariantly
finite subcategory of ${\mathscr A}$. We claim that $\varepsilon_{B}\circ F(f)$: $F(C)\longrightarrow FG(B)\longrightarrow B$ is a right $F(\mathscr C)$-approximation,
i.e., $F({\mathscr C})$ is a contravariantly finite subcategory of ${\mathscr B}$. In fact, we only need to find a map that satisfies the following
diagram is commutative:
$$\xymatrix{&&F(X)\ar[d]^{g}\ar@{.>}[dll]\\
F(C)\ar[rr]^{\varepsilon_{B}\circ F(f)}&&B
}$$
for any $g$ and $X\in\mathscr C$. Since $f$ is a right $\mathscr C$-approximation, we have the following commutative diagram:
$$\xymatrix{
X\ar[r]^{\eta_{X}}\ar@{.>}[d]^{h}&GF(X)\ar[d]^{G(g)}\\
C\ar[r]^{f}&G(B)
}$$
i.e., $fh=G(g)\eta_{X}$. We consider the following commutative diagram:
$$\xymatrix{
FGF(X)\ar[r]^{\varepsilon_{F(X)}}\ar[d]^{FG(g)}&F(X)\ar[d]^{g}\\
FG(B)\ar[r]^{\varepsilon_{B}}&B
}$$
i.e., $g\varepsilon_{F(X)}=\varepsilon_{B}FG(g)$. By the Lemma \ref{adjoindpro}, $g=g\varepsilon_{F(X)}F(\eta_{X})=\varepsilon_{B}FG(g)F(\eta_{X})
=\varepsilon_{B}F(G(g)\eta_{X})=\varepsilon_{B}F(fh)=\varepsilon_{B}F(f)F(h)$.
So there is a morphism $F(h)$: $F(X)\longrightarrow F(C)$ such that $g=\varepsilon_{B}\circ F(f)F(h)$, i.e.,
$F({\mathscr C})$ is a contravariantly finite subcategory of ${\mathscr B}$.

Note that $G$ is faithful if and only if the counit $\varepsilon_{B}$: $FG(B)\longrightarrow B$ is surjective for any $B \in \mathscr B$.
As ${\mathscr C}$ is admissible, $f$ is an epimorphism. By the discuss above, $\varepsilon_{B}F(f)$ is surjective since $F$ preserves epimorphism.
i.e., $F({\mathscr C})$ is also admissible.
\ \hfill $\Box$

\mskip\

%
%

We don't know whether the converse of Lemma \ref{YL1} is true or not, but we have the following conclusion.

\bg{Pro}\label{YL2}
Let $F:$ $\mathscr A \longrightarrow \mathscr B$ be a functor, ${\mathscr C}$ be a subcategory of ${\mathscr A}$ such that $F(\mathscr C)$ is a
covariantly $($resp., contravariantly$)$ finite subcategory of ${\mathscr B}$. If F is fully faithful, then ${\mathscr C}$ is
a covariantly $($resp., contravariantly$)$ finite subcategory of ${\mathscr A}$.
\ed{Pro}

\Pf. For any $A\in\mathscr A$, there is a left $F(\mathscr C)$-approximation $\bar{f}$: $F(A)\longrightarrow F(C)$ since $F(\mathscr C)$ is
covariantly finite. Note that $\mathrm{Hom}_{\mathscr B}(F(A),~F(C))\cong \mathrm{Hom}_{\mathscr A}(A,~C)$ since $F$ is fully faithful. We claim that
$f$: $A\longrightarrow C$ is a left $\mathscr C$-approximation, which is image of $\bar{f}$. As $\mathrm{Hom}_{\mathscr B}(\bar{f},~F(C'))$ is surjective
for any $C'\in\mathscr C$, $\mathrm{Hom}_{\mathscr A}(f,~C')$ is also surjective, i.e., ${\mathscr C}$ is a covariantly finite subcategory of ${\mathscr A}$.
\ \hfill $\Box$

\mskip\
In order to prove the main results of this section, we also need the following conclusion.

\bg{Pro}\label{pro1}
Let $\mathscr A$ and $\mathscr B$ be abelian categories, and let $F$: $\mathscr B\longrightarrow\mathscr A$, $G$: $\mathscr A\longrightarrow\mathscr B$,
$H$: $\mathscr B\longrightarrow\mathscr A$ be additive functors such that ($F$, $G$, $H$) is an adjoint triple. If ($\mathscr U$, $\mathscr V$) is a balanced
pair in $\mathscr B$, then (F($\mathscr U$), H($\mathscr V$)) is a balanced pair in $\mathscr A$. If, in addition, $G$ is faithful and ($\mathscr U$, $\mathscr V$)
is admissible in $\mathscr B$, then (F($\mathscr U$), H($\mathscr V$)) is admissible in $\mathscr A$.
\ed{Pro}

\Pf. By the lemma \ref{YL1}, $F({\mathscr U})$ is contravariantly finite subcategory of ${\mathscr A}$.
Let $M^{\bullet}$ be a complex in $\mathscr A$. $\mathrm{Hom}_{\mathscr A}(F({\mathscr U}),~M^{\bullet})$ is acyclic $\Longleftrightarrow$
$\mathrm{Hom}_{\mathscr A}({\mathscr U},~G(M^{\bullet}))$ is acyclic (since $(F,~G)$ is a adjoint pair) $\Longleftrightarrow$
$\mathrm{Hom}_{\mathscr A}(G(M^{\bullet}),~{\mathscr V})$ is acyclic (since $({\mathscr U}$, ${\mathscr V})$ is a balanced pair) $\Longleftrightarrow$
$\mathrm{Hom}_{\mathscr A}(M^{\bullet},~H({\mathscr V}))$ is acyclic (since $(G,~H)$ is a adjoint pair).
So the pair $(F(\mathscr U$), $H(\mathscr V$)) is balanced by the proposition 2.2 in \cite{CHEN}.

In addition, if $G$ is faithful and ($\mathscr U$, $\mathscr V$) is admissible, by the lemma \ref{YL1},
we have that $F(\mathscr U)$ is admissible. i.e., (F($\mathscr U$), H($\mathscr V$)) is admissible in $\mathscr A$.
\ \hfill $\Box$

\bg{Th}\label{result1}
Suppose that $({\mathscr X'}$, ${\mathscr Y'})$ and $({\mathscr X''}$, ${\mathscr Y''})$ are two balanced pairs in ${\mathscr A'}$ and ${\mathscr A''}$,
respectively, and let
$${\mathscr X}=\{X\in{\mathscr A}|i^{\ast}X\in {\mathscr X'},~j^{\ast}X\in {\mathscr X''}\}$$
$${\mathscr Y}=\{Y\in{\mathscr A}|i^{!}Y\in {\mathscr Y'},~j^{\ast}Y\in {\mathscr Y''}\}$$
then the following statements hold:

$(1)$ $i^{\ast}{\mathscr X}={\mathscr X'}$, $i^{!}{\mathscr Y}={\mathscr Y'}$, $j^{\ast}{\mathscr X}={\mathscr X''}$ and $j^{\ast}{\mathscr Y}={\mathscr Y''}$;

$(2)$ the pair $(j_{!}j^{\ast}{\mathscr X}$, $j_{\ast}j^{\ast}{\mathscr Y})$ is balanced pair.
Moreover, if $({\mathscr X''}$, ${\mathscr Y''})$ is admissible and $j^{\ast}$ is faithful, then $(j_{!}j^{\ast}{\mathscr X}$, $j_{\ast}j^{\ast}{\mathscr Y})$ is also admissible;

$(3)$ if $({\mathscr X''}$, ${\mathscr Y''})$ is of finite dimensions and $j^{\ast}$ is faithful,
then the balanced pair $(j_{!}j^{\ast}{\mathscr X}$, $j_{\ast}j^{\ast}{\mathscr Y})$ is of finite dimensions.

\ed{Th}

\Pf. $(1)$ Obviously, $i^{\ast}{\mathscr X}\subseteq{\mathscr X'}$. By the proposition \ref{propR}, we have that $\mathscr X'\cong i^{\ast}i_{\ast}\mathscr X'$
and $i_{\ast}\mathscr X'\subseteq\mathscr X$, i.e., $\mathscr X'\subseteq i^{\ast} \mathscr X$.
So $\mathscr X'= i^{\ast}i_{\ast}\mathscr X'$. The other three equations are similar.

$(2)$ Note that $(j_{!}j^{\ast}{\mathscr X},~j_{\ast}j^{\ast}{\mathscr Y})=(j_{!}{\mathscr X''},~j_{\ast}{\mathscr Y''})$. According to the proposition
\ref{pro1}, we can easily prove that it is correct.

$(3)$ Set $\mathscr X''$-dim $\mathscr A''$=$n$. For any $A\in\mathscr A$, since $j_{!}j^{\ast}{\mathscr X}$ is contravariantly finite subcategory,
there is a $\mathrm{Hom}(j_{!}j^{\ast}{\mathscr X},~-)$-cyclic complex
$$\xymatrix{
\cdots\ar[r]&j_{!}j^{\ast}X_{n}\ar[r]^{d_{n}}&\cdots\ar[r]&j_{!}j^{\ast}X_{1}\ar[r]&j_{!}j^{\ast}X_{0}\ar[r]&A\ar[r]&0
}$$
where each $X_{i}\in \mathscr X$. i.e., the sequence
$$\xymatrix@C=0.5cm{
\cdots\ar[r]&\mathrm{Hom}(j_{!}j^{\ast}{\mathscr X},~j_{!}j^{\ast}X_{n})\ar[r]&\cdots\ar[r]
&\mathrm{Hom}(j_{!}j^{\ast}{\mathscr X},~j_{!}j^{\ast}X_{0})\ar[r]&\mathrm{Hom}(j_{!}j^{\ast}{\mathscr X},~A)\ar[r]&0
}$$
is exact. Applying the exact functor $j^{\ast}$ to above exact sequence, we can obtain the following exact sequence
$$\xymatrix{
\cdots\ar[r]&\mathrm{Hom}(j^{\ast}{\mathscr X},~j^{\ast}X_{n})\ar[r]&\cdots\ar[r]
&\mathrm{Hom}(j^{\ast}{\mathscr X},~j^{\ast}X_{0})\ar[r]&\mathrm{Hom}(j^{\ast}{\mathscr X},~j^{\ast}A)\ar[r]&0
}$$
since $\mathrm{Id}_{\mathscr A''}\cong j^{\ast}j_{!}$ by Proposition \ref{propR}. Consequently, the following exact sequence gives a
left $j^{\ast}{\mathscr X}=\mathscr X''$-resolution of $j^{\ast}A$.
$$\xymatrix{
\cdots\ar[r]&j^{\ast}X_{n}\ar[r]^{j^{\ast}d_{n}}&\cdots\ar[r]&j^{\ast}X_{1}\ar[r]
&j^{\ast}X_{0}\ar[r]&j^{\ast}A\ar[r]&0
}$$
By the lemma 2.4 in \cite{CHEN}, $j^{\ast}\mathrm{Ker}(d_{n-1})=\mathrm{Ker}(j^{\ast}d_{n-1})\in\mathscr X''=j^{\ast}{\mathscr X}=j^{\ast}j_{!}j^{\ast}{\mathscr X}$ since $\mathscr X''$-dim $\mathscr A''$=$n$.
Note that the functor $j^{\ast}$ is exact and faithful. Thus $\mathrm{Ker}(d_{n-1})\in j_{!}j^{\ast}{\mathscr X}$.
From the first sequence, we can know that $j_{!}j^{\ast}{\mathscr X}$-dim $A$ is finite, i.e.,
the pair $(j_{!}j^{\ast}{\mathscr X}$, $j_{\ast}j^{\ast}{\mathscr Y})$ is of finite dimensions.
\ \hfill $\Box$

\mskip\

It is very interesting to prove that the pair $({\mathscr X}$, ${\mathscr Y})$ defined in Theorem \ref{result1} is a balanced
pair. Unfortunately, we don't know wether this conclusion is correct or not.
%
%
%
%

\bg{Lem}\label{YL3}
Let $\mathscr A$ and $\mathscr B$ be abelian categories, and let
$$L^{\bullet}=:\xymatrix{\cdots\ar[r]&L_{-1}\ar[r]^{d_{-1}}&L_{0}\ar[r]^{d_{0}}&L_{1}\ar[r]&\cdots}$$
be in${\mathscr A}$. If the functor $F:$ $\mathscr A \longrightarrow \mathscr B$ be exact and faithful. Then $L^{\bullet}$ is acyclic if and only if $F(L^{\bullet})$ is acyclic.
\ed{Lem}

\Pf. The necessity is obvious since $F$ is exact.

Since $F(L^{\bullet})=:$ $\xymatrix{\cdots\ar[r]&FL_{-1}\ar[r]^{Fd_{-1}}&FL_{0}\ar[r]^{Fd_{0}}&FL_{1}\ar[r]&\cdots}$ is acyclic and $F$ is exact,
then the $i$-th homology  0=$\mathrm{H}_{i}(F(L^{\bullet}))=\mathrm{Ker} Fd_{i}/\mathrm{Im} Fd_{i-1}\cong F(\mathrm{Ker} d_{i}/\mathrm{Im} d_{i-1})$.
Note that an exact functor $F$ is faithful if and only if  it does not take any non-zero object to zero. So $\mathrm{Ker} d_{i}/\mathrm{Im} d_{i-1}=0$.
i.e., $L^{\bullet}$ is acyclic.
\ \hfill $\Box$

\bg{Th}\label{result2}
Suppose that $({\mathscr X}$, ${\mathscr Y})$ is a balanced pair in ${\mathscr A}$ and $j^{\ast}$ is faithful, then we have

$(1)$ the two pairs $(i^{\ast}{\mathscr X}$, $i^{!}{\mathscr Y})$
and $(j^{\ast}{\mathscr X}$, $j^{\ast}{\mathscr Y})$ are balanced in ${\mathscr A'}$ and ${\mathscr A''}$,
respectively;

$(2)$ if $({\mathscr X}$, ${\mathscr Y})$ is admissible, then the two pairs $(i^{\ast}{\mathscr X}$, $i^{!}{\mathscr Y})$
and $(j^{\ast}{\mathscr X}$, $j^{\ast}{\mathscr Y})$ are admissible;

$(3)$ if $({\mathscr X}$, ${\mathscr Y})$ is of finite dimensions,
then the balanced pair $(j^{\ast}{\mathscr X}$, $j^{\ast}{\mathscr Y})$ is of finite dimensions.
In addition, if ${\mathscr X}=i_{\ast}i^{\ast}\mathscr X$, then $(i^{\ast}{\mathscr X}$, $i^{!}{\mathscr Y})$ is also of finite dimensions.
\ed{Th}

\Pf. $(1)$ The first statement holds by the proposition \ref{pro1}.
By the lemma \ref{YL1}, $j^{\ast}{\mathscr X}$ is a contravariantly finite subcategories.


Let $N^{\bullet}$ be a complex in $\mathscr A''$. $\mathrm{Hom}_{\mathscr A''}(j^{\ast}{\mathscr X},~N^{\bullet})$ is acyclic $\Longleftrightarrow$
$\mathrm{Hom}_{\mathscr A''}({\mathscr X},~j_{\ast}N^{\bullet})$ is acyclic (since $(j^{\ast},~j_{\ast})$ is a adjoint pair) $\Longleftrightarrow$
$\mathrm{Hom}_{\mathscr A''}(j_{\ast}N^{\bullet},~{\mathscr Y})$ is acyclic (since $({\mathscr X}$, ${\mathscr Y})$ is a balanced pair) $\Longleftrightarrow$
$\mathrm{Hom}_{\mathscr A''}(j^{\ast}j_{\ast}N^{\bullet},~j^{\ast}{\mathscr Y})$ is acyclic (since $j^{\ast}$ is faithful, exact and the lemma \ref{YL3}) $\Longleftrightarrow$
$\mathrm{Hom}_{\mathscr A''}(N^{\bullet},~j^{\ast}{\mathscr Y})$ is acyclic (since $j^{\ast}j_{\ast}\cong \mathrm{Id}\mathscr A''$).
So the pair $(j^{\ast}{\mathscr X}$, $j^{\ast}{\mathscr Y})$ is balanced.

$(2)$ The pair $(i^{\ast}{\mathscr X}$, $i^{!}{\mathscr Y})$ is admissible by the proposition \ref{pro1}.

Let $Q^{\bullet}$ be a complex in $\mathscr A''$. $\mathrm{Hom}_{\mathscr A''}(j^{\ast}{\mathscr X},~Q^{\bullet})$ is acyclic $\Longrightarrow$
$\mathrm{Hom}_{\mathscr A''}({\mathscr X},~j_{\ast}Q^{\bullet})$ is acyclic (since $(j^{\ast},~j_{\ast})$ is a adjoint pair) $\Longrightarrow$
$j_{\ast}Q^{\bullet}$ is acyclic (since $({\mathscr X}$, ${\mathscr Y})$ is admissible) $\Longrightarrow$
$j^{\ast}j_{\ast}Q^{\bullet}$ is acyclic (since $j^{\ast}$ is exact) $\Longrightarrow$
$Q^{\bullet}$ is acyclic (since $j^{\ast}j_{\ast}\cong \mathrm{Id}\mathscr A''$).
So the pair $(j^{\ast}{\mathscr X}$, $j^{\ast}{\mathscr Y})$ is admissible.

$(3)$ For any $A'' \in {\mathscr A''}$, there is a left ${\mathscr X}$-resolution of $j_! A''$ since $({\mathscr X}$, ${\mathscr Y})$ is of finite dimensions.
$$\xymatrix{
0\ar[r]&X_{n}\ar[r]&\cdots\ar[r]&X_{1}\ar[r]&X_{0}\ar[r]&j_!A''\ar[r]&0
}$$
with $X_{i}\in{\mathscr X}$ for any $i$. i.e., the following sequence is acyclic.
$$\xymatrix{
0\ar[r]&\mathrm{Hom}_{\mathscr A}(\mathscr X,~X_{n})\ar[r]&\cdots\ar[r]
&\mathrm{Hom}_{\mathscr A}(\mathscr X,~X_{0})\ar[r]&\mathrm{Hom}_{\mathscr A}(\mathscr X,~j_!A'')\ar[r]&0
}$$
Since $j^{\ast}$ is exact, then the sequence
$$\xymatrix@C=0.5cm{
0\ar[r]&\mathrm{Hom}_{\mathscr A''}(j^{\ast}\mathscr X,~j^{\ast}X_{n})\ar[r]&\cdots\ar[r]
&\mathrm{Hom}_{\mathscr A''}(j^{\ast}\mathscr X,~j^{\ast}X_{0})\ar[r]&
\mathrm{Hom}_{\mathscr A''}(j^{\ast}\mathscr X,~j^{\ast}j_!A'')\ar[r]&0
}$$
is also acyclic.
Note that $\mathrm{Hom}_{\mathscr A''}(j^{\ast}\mathscr X,~j^{\ast}j_!A'')\cong\mathrm{Hom}_{\mathscr A''}(j^{\ast}\mathscr X,~A'')$
by the proposition \ref{propR}. It is easy to see that the sequence
$$\xymatrix{
0\ar[r]&j^{\ast}X_{n}\ar[r] &\cdots\ar[r]&j^{\ast}X_{1}\ar[r]&j^{\ast}X_{0}\ar[r]& A''\ar[r]&0
}$$
is a left $j^{\ast}{\mathscr X}$-resolution of $ A''$. So the balanced pair $(j^{\ast}{\mathscr X}$, $j^{\ast}{\mathscr Y})$ is of finite dimensions.

For any $A'\in {\mathscr A'}$, since $i^{\ast}{\mathscr X}$ is contravariantly finite subcategory, there is a left $i^{\ast}{\mathscr X}$-resolution of $A'$.
$$\xymatrix{
\cdots\ar[r]&i^{\ast}X_{n}\ar[r]^{\partial_{n}}&\cdots\ar[r]&i^{\ast}X_{1}\ar[r]&i^{\ast}X_{0}\ar[r]&A'\ar[r]&0
}$$
with $X_{i}\in{\mathscr X}$ for any $i$. i.e., the following sequence is acyclic.
$$\xymatrix{
\cdots\ar[r]&\mathrm{Hom}_{\mathscr A'}(i^{\ast}\mathscr X,~i^{\ast}X_{n})\ar[r]&\cdots\ar[r]
&\mathrm{Hom}_{\mathscr A'}(i^{\ast}\mathscr X,~i^{\ast}X_{0})\ar[r]&\mathrm{Hom}_{\mathscr A'}(i^{\ast}\mathscr X,~A')\ar[r]&0
}$$
Since $(i^{\ast},~i_{\ast})$ is adjoint pair, then the sequence
$$\xymatrix{
\cdots\ar[r]&\mathrm{Hom}_{\mathscr A}(\mathscr X,~i_{\ast}i^{\ast}X_{n})\ar[r]&\cdots\ar[r]
&\mathrm{Hom}_{\mathscr A}(\mathscr X,~i_{\ast}i^{\ast}X_{0})\ar[r]&\mathrm{Hom}_{\mathscr A}(\mathscr X,~i_{\ast}A')\ar[r]&0
}$$
is acyclic. i.e., the sequence
$$\xymatrix{
\cdots\ar[r]&i_{\ast}i^{\ast}X_{n}\ar[r]^{i_{\ast}\partial_{n}}&\cdots\ar[r]
&i_{\ast}i^{\ast}X_{0}\ar[r]&i_{\ast}A'\ar[r]&0
}$$
is a left ${\mathscr X}$-resolution of $i_{\ast}A'$. Since $({\mathscr X}$, ${\mathscr Y})$ is of finite dimensions, by the lemma 2.4 in \cite{CHEN},
$i_{\ast} (\mathrm{Ker}\partial_{n-1})\cong\mathrm{Ker} (i_{\ast}\partial_{n-1})\in \mathscr X=i_{\ast}i^{\ast}\mathscr X$.
Note that the functor $i_{\ast}$ are faithful and exact. Thus we have that $\mathrm{Ker}\partial_{n}\in i^{\ast}\mathscr X$,
since the exact and faithful functor does not take any non-zero object to zero.
So $i^{\ast}{\mathscr X}$-dim $A' \leq n$. i.e., $(i^{\ast}{\mathscr X}$, $i^{!}{\mathscr Y})$ is also of finite dimensions.
\ \hfill $\Box$

\mskip\

In fact, from the proof above, we have that $j^{\ast}{\mathscr X}$-dim $\mathscr A''\leq$ ${\mathscr X}$-dim $\mathscr A$
and $i^{\ast}{\mathscr X}$-dim $\mathscr A'\leq$ ${\mathscr X}$-dim $\mathscr A$.

\section{Relative tilting modules}
\vskip 10pt

In this section, we also assume that $({\mathscr X}$, ${\mathscr Y})$ is an admissible balanced pair in the Abelian category ${\mathscr A}$.
Whenever the subcategory ${\mathscr X}$ (respectively, ${\mathscr Y}$) appears, we always suppose
that it is a part of the admissible balanced pair $({\mathscr X}$, ${\mathscr Y})$.
Here, we will give the notion of $n$-${\mathscr X}$-tilting module with respect to $({\mathscr X}$, ${\mathscr Y})$,
and obtain a characterization of $n$-${\mathscr X}$-tilting module, which similar to Bazzoni characterization of $n$-tilting module \cite{BS}.

Note that the functor $\mathrm{Ext}^{i}_{{\mathscr A}}(-,~-)$ is based on the classical balanced pair $(\mathrm{Proj} R,~\mathrm{Inj} R)$, it induces
an isomorphism of cohomology groups whether we take a projective resolution of the first variable or take
an injective coresolution of the second variable. From this viewpoint, we have the functor $\mathrm{Ext}^{i}_{{\mathscr X}}(-,~-)$ with respect
to admissible balanced pair $({\mathscr X}$, ${\mathscr Y})$. By Lemma 2.1 \cite{CHEN}, the functor is well-define.
Firstly, we give the following some notions which are widely used in this section.

$\mathrm{Add}T$ =: $\{K|$ there exists an object $J$ such that $K\bigoplus J\cong T^{(I)}$ for some set $I. \}$

$\widehat{\mathrm{Add}_{\mathscr X}T}$ =: $\{L|$ for some $n$, there is $\mathscr X$-exact sequence $0\longrightarrow T_{n}\longrightarrow \cdots\longrightarrow 
T_{1}\longrightarrow L\longrightarrow0$ with $T_{i}\in\mathrm{Add}T$ for all $1\leq i\leq n.\}$

$(\widehat{\mathrm{Add}_{\mathscr X}T})_{n}$ =: $\{L|$ there is $\mathscr X$-exact sequence $0\longrightarrow T_{n}\longrightarrow \cdots\longrightarrow
T_{1}\longrightarrow L\longrightarrow0$ with $T_{i}\in\mathrm{Add}T$ for all $1\leq i\leq n.\}$

$\mathrm{Pres}^{n}_{\mathscr X}T$ =: $\{L|$ there is $\mathscr X$-exact sequence $T_{n}\longrightarrow \cdots\longrightarrow T_{1}\longrightarrow L\longrightarrow0$
with $T_{i}\in\mathrm{Add}T$ for all $1\leq i\leq n.\}$

$T^{\mathscr X \perp}$=: $\{M|$ $\mathrm{Ext}^{i\geq1}_{\mathscr X}(T,~M)=0.\}$

$_{T}\mathscr X$=: $\{M|$ there is $\mathscr X$-exact sequence $\xymatrix{\cdots \ar[r] &T_{n}\ar[r]^{f_{n}}&\cdots\ar[r]&T_{1}\ar[r]^{f_{1}}&T_{0}\ar[r]^{f_{0}}&M\ar[r]&0}$
with $T_{i}\in\mathrm{Add}T$ and $\mathrm{Im}f_{i}\in T^{\mathscr X \perp}$ for all $i$ $.\}$

Dually, we can define these notions: $\widecheck{\mathrm{Add}_{\mathscr X}T}$, $(\widecheck{\mathrm{Add}_{\mathscr X}T})_{n}$, $^{\mathscr X \perp}T$ and $\mathscr X _{T}$.

\bg{Def}\label{tilting}%
Let ${\mathscr A}$ be an Abelian category and $T\in {\mathscr A}$. Then $T$ is called $n$-${\mathscr X}$ tilting (with respect to $({\mathscr X}$, ${\mathscr Y})$) if it satisfying the following conditions:

$(1)$ ${\mathscr X}$-dim $T \leq n$;

$(2)$ $T$ is ${\mathscr X}$-self-orthogonal, i.e., $\mathrm{Ext}^{i}_{{\mathscr X}}(T,~T^{(I)})=0$, for each $i>0$ and all sets $I$;

$(3)$ there is a ${\mathscr X}$-exact sequence (i.e., the sequence is ${\mathscr X}$-acyclic)
$$\xymatrix{
0\ar[r]&X\ar[r]&T_{0}\ar[r]&\cdots\ar[r]&T_{n}\ar[r]&0
}$$
for any $X\in {\mathscr X}$, where $T_{i}\in \mathrm{Add}T$ for any $i$.
\ed{Def}

If $T$ satisfies the above conditions $(1)$ and $(2)$, then $T$ is said to be partial $n$-${\mathscr X}$ tilting.

\bg{Exm}
$(1)$ It is easily verify that $(\mathrm{Proj} R,~\mathrm{Inj} R)$ is an admissible balanced pair in $\mathrm{Mod}$-R.
Then the n-tilting module \cite{WEI} is a n-$\mathrm{Proj}$ tilting module with respect to balanced pair $(\mathrm{Proj} R,~\mathrm{Inj} R)$.

$(2)$ By Proposition 2.6 in \cite{CHEN}, we can obtain that $(\mathcal{GP},~\mathcal{GI})$ in Example \ref{eg1} is an admissible balanced pair. Then the $n$-Gorenstein tilting modules \cite{YLO} is
n-$\mathcal{GP}$ tilting module with respect to balanced pair $(\mathcal{GP},~\mathcal{GI})$.
\ed{Exm}

\bg{Lem}\label{laba}
Let ${\mathscr A}$ be an Abelian category and $T\in {\mathscr A}$.

$(1)$ If $U\in \mathrm{Gen}_{\mathscr X} T$, then there is an $\mathscr X$-exact sequence
$\xymatrix{0\ar[r]&V\ar[r]&T_{U}\ar[r]&U\ar[r]&0}$ with $T_{U}\in \mathrm{Add}T$,
which stays exactness after applying the functor $\mathrm{Hom}(T,~-)$.

$(2)$ If $T$ is ${\mathscr X}$-self-orthogonal, and $T^{\mathscr X\bot}\subseteq \mathrm{Gen}_{\mathscr X} T (=\mathrm{Pres}^{1}_{\mathscr X} T)$,
then $T^{\mathscr X\bot}=$ $_{T}\mathscr X$.

$(3)$ If T is $n$-${\mathscr X}$-tilting, then $T^{\mathscr X\bot}\subseteq \mathrm{Gen}_{\mathscr X} T$.
Specially, $T^{\mathscr X\bot}=$ $_{T}\mathscr X$.
\ed{Lem}

\Pf. (1) There is an $\mathscr X$-exact sequence
$\xymatrix{0\ar[r]&V'\ar[r]&T'_{U}\ar[r]&U\ar[r]&0}$ with $T_{U}\in \mathrm{Add}T$ since $U\in \mathrm{Gen}_{\mathscr X} T$.
So the evaluation mapping $f:$ $T_{U}\longrightarrow U$ is surjective. We can obtain the following commutative diagram
$$\xymatrix{
0\ar[r]&V'\ar[r]\ar[d]&T'_{U}\ar[d]\ar[r]&U\ar@^{=}[d]\ar[r]&0\\
0\ar[r]&V\ar[r]&T_{U}\ar[r]^{f}&U\ar[r]&0
}$$
It is easy to see that the diagram is a pushout. Note that every pushout and pullback of an $\mathscr X$-exact sequence is again $\mathscr X$-exact
by \cite[section 1]{AS}. Consequence, the second row in above diagram is our desired.

(2) Clearly, $_{T}\mathscr X\subseteq T^{\mathscr X\bot}$. For any $U\in T^{\mathscr X\bot}\subseteq \mathrm{Gen}_{\mathscr X} T$,
there is an $\mathscr X$-exact sequence $\xymatrix{0\ar[r]&V\ar[r]&T_{U}\ar[r]&U\ar[r]&0}$ with $T_{U}\in \mathrm{Add}T$,
and it stays exactness after applying the functor $\mathrm{Hom}(T,~-)$ by (1). It is not difficult to prove that
$V\in T^{\mathscr X\bot}$. Repeating this process to $V$, we can prove that $U\in ~_{T}\mathscr X$. i.e.,
$T^{\mathscr X\bot}=$ $_{T}\mathscr X$.

(3) For any $U\in T^{\mathscr X\bot}$, we take an $\mathscr X$-exact sequence
$\xymatrix{0\ar[r]&U_{1}\ar[r]&X_{U}\ar[r]&U\ar[r]&0}$ with $X_{U}\in \mathscr X$ since $\mathscr X$ is
contravariantly finite and admissible. Since $T$ is $n$-${\mathscr X}$-tilting, there is an $\mathscr X$-exact sequence
$\xymatrix{0\ar[r]&X_{U}\ar[r]&T_{0}\ar[r]&X'_{U}\ar[r]&0}$ with $T_{0}\in \mathrm{Add}T$ and $X'_{U}\in \widecheck{\mathrm{Add}_{\mathscr X}T}$.
Therefore we have the following commutative diagram.
$$\xymatrix{
&&0\ar[d]&0\ar[d]\\
0\ar[r]&U_{1}\ar[r]\ar@^{=}[d]&X_{U}\ar[d]\ar[r]&U\ar[d]\ar[r]&0\\
0\ar[r]&U_{1}\ar[r]&T_{0}\ar[d]\ar[r]&V\ar[r]\ar[d]&0\\
&&X'_{U}\ar[d]\ar@^{=}[r]&X'_{U}\ar[d]\\
&&0&0
}$$
It is easy to verify that the third column is split. i.e., $V=U\bigoplus X'_{U}$. Note that all sequence in above are $\mathscr X$-exact.
From the second row, we have that $V\in\mathrm{Gen}_{\mathscr X} T$, i.e., $U\in\mathrm{Gen}_{\mathscr X} T$.
Consequence, $T^{\mathscr X\bot}\subseteq \mathrm{Gen}_{\mathscr X} T$. Specially, $T^{\mathscr X\bot}=$ $_{T}\mathscr X$ by (2).
\ \hfill $\Box$

\bg{Lem}\label{ltbx}
Let $T$ be ${\mathscr X}$-self-orthogonal. Then

$(1)$ $_{T}\mathscr X$ is closed under ${\mathscr X}$-extensions and closed under direct summands.

$(2)$ $\widehat{\mathrm{Add}_{\mathscr X}T}=$ $_{T}\mathscr X\bigcap T^{{\mathscr X}\bot}$. Specially,
$\widehat{\mathrm{Add}_{\mathscr X}T}$ is closed under direct summands.
\ed{Lem}

\Pf. (1) The prove of first statement similar to the usual horseshoe lemma. Now suppose that
$U=M\bigoplus N$ is in $_{T}\mathscr X$. There is an ${\mathscr X}$-exact sequence
$0\longrightarrow U'\longrightarrow T_{U}\longrightarrow U\longrightarrow 0$ with $U'\in$ $_{T}\mathscr X$
and $T_{U}\in\mathrm{Add}T$. We consider the following commutative diagram.
$$\xymatrix{
&0\ar[d]&0\ar[d]\\
&U'\ar@^{=}[r]\ar[d]&U'\ar[d]\\
0\ar[r]&X\ar[r]\ar[d]&T_{U}\ar[d]\ar[r]&N\ar@^{=}[d]\ar[r]&0\\
0\ar[r]&M\ar[d]\ar[r]&U\ar[d]\ar[r]&N\ar[r]&0\\
&0&0
}$$
We claim that $X\bigoplus N$ $\in _{T}\mathscr X$. Indeed, the sequence
$0\longrightarrow U'\longrightarrow X\bigoplus N\longrightarrow M\bigoplus N=U\longrightarrow 0$
is ${\mathscr X}$-exact. By the first statement, $X\bigoplus N$ $\in _{T}\mathscr X$.
Hence we can deduce by recursiveness that $_{T}\mathscr X$ is closed under direct summands.

(2) The proof of this equation is trivial.
\ \hfill $\Box$

\bg{Lem}\label{ltba}
Let $T$ be ${\mathscr X}$-self-orthogonal. The sequence $0\longrightarrow V\longrightarrow M_{n}\longrightarrow \cdots\longrightarrow M_{1}\longrightarrow U\longrightarrow 0$ is ${\mathscr X}$-exact, where $M_{i}\in$ $_{T}{\mathscr X}$. Then

$(1)$ there is an ${\mathscr X}$-exact sequence $0\longrightarrow A_{n}\longrightarrow B_{n}\longrightarrow V\longrightarrow 0$, where
$A_{n}\in$ $_{T}{\mathscr X}$ and $B_{n}$ satisfying there is an ${\mathscr X}$-exact sequence
$0\longrightarrow B_{n}\longrightarrow T_{n}\longrightarrow \cdots\longrightarrow T_{1}\longrightarrow U\longrightarrow 0$
with $T_{i}\in \mathrm{Add}T$.

$(2)$ If $U\in$ $_{T}{\mathscr X}$, then there is an ${\mathscr X}$-exact sequence
$0\longrightarrow A\longrightarrow B\longrightarrow V\longrightarrow 0$ with $A\in$ $_{T}{\mathscr X}$
and $B\in(\widecheck{\mathrm{Add}_{\mathscr X}T})_{n}$.

\ed{Lem}

\Pf. To prove (1), we use induction to $n$. When $n=1$. Note that there is an ${\mathscr X}$-exact sequence
$0\longrightarrow A_{1}\longrightarrow T_{1}\longrightarrow M_{1}\longrightarrow 0$ with $A_{1}\in$ $_{T}{\mathscr X}$
and $T_{1}\in\mathrm{Add}T$ since $M_{1}\in$ $_{T}{\mathscr X}$. We consider the following commutative diagram.
$$\xymatrix{
&0\ar[d]&0\ar[d]\\
&A_{1}\ar@^{=}[r]\ar[d]&A_{1}\ar[d]\\
0\ar[r]&B_{1}\ar[r]\ar[d]&T_{1}\ar[d]\ar[r]&U\ar@^{=}[d]\ar[r]&0\\
0\ar[r]&V\ar[d]\ar[r]&M_{1}\ar[d]\ar[r]&U\ar[r]&0\\
&0&0
}$$
Then the first column and the second row is our desired.

Assume that the statement holds for $n-1$. Set $\mathrm{Coker}(V\longrightarrow M_{n})=V'$. By assumption, we have an
${\mathscr X}$-exact sequence $0\longrightarrow A'_{n-1}\longrightarrow B'_{n-1}\longrightarrow V'\longrightarrow 0$, where
$A'_{n-1}\in$ $_{T}{\mathscr X}$. We consider the following commutative diagram.
$$\xymatrix{
&&0\ar[d]&0\ar[d]\\
&&A'_{n-1}\ar@^{=}[r]\ar[d]&A'_{n-1}\ar[d]\\
0\ar[r]&V\ar[r]\ar@^{=}[d]&X\ar[d]\ar[r]&B'_{n-1}\ar[d]\ar[r]&0\\
0\ar[r]&V\ar[r]&M_{n}\ar[d]\ar[r]&V'\ar[d]\ar[r]&0\\
&&0&0
}$$
It is easy to see that all sequence in above are ${\mathscr X}$-exact. By Lemma \ref{ltbx}, $X\in $ $_{T}{\mathscr X}$.
Thus there is an ${\mathscr X}$-exact sequence $0\longrightarrow X\longrightarrow T_{n}\longrightarrow A_{n}\longrightarrow 0$, where
$A_{n}\in$ $_{T}{\mathscr X}$ and $T_{X}\in \mathrm{Add}T$. We consider the following commutative diagram.
$$\xymatrix{
&0\ar[d]&0\ar[d]\\
&A_{n}\ar@^{=}[r]\ar[d]&A_{n}\ar[d]\\
0\ar[r]&B_{n}\ar[r]\ar[d]&T_{n}\ar[d]\ar[r]&B'_{n-1}\ar@^{=}[d]\ar[r]&0\\
0\ar[r]&V\ar[d]\ar[r]&X\ar[d]\ar[r]&B'_{n-1}\ar[r]&0\\
&0&0
}$$
Then the first column and the second row is what we want.

(2) Since $U\in$ $_{T}{\mathscr X}$, then there is an ${\mathscr X}$-exact sequence
$0\longrightarrow U_{1}\longrightarrow T_{U}\longrightarrow U\longrightarrow 0$, where
$U_{1}\in$ $_{T}{\mathscr X}$ and $T_{U}\in \mathrm{Add}T$. Set $\mathrm{Ker}(M_{1}\longrightarrow U)=V_{1}$. We consider the following commutative diagram.
$$\xymatrix{
&&0\ar[d]&0\ar[d]\\
&&U_{1}\ar@^{=}[r]\ar[d]&U_{1}\ar[d]\\
0\ar[r]&V_{1}\ar[r]\ar@^{=}[d]&N_{1}\ar[d]\ar[r]&T_{U}\ar[d]\ar[r]&0\\
0\ar[r]&V_{1}\ar[r]&M_{1}\ar[d]\ar[r]&U\ar[d]\ar[r]&0\\
&&0&0
}$$
It is easy to see that all sequence in above are ${\mathscr X}$-exact. By Lemma \ref{ltbx}, $N_{1}\in $ $_{T}{\mathscr X}$.
Thus we can obtain an ${\mathscr X}$-exact sequence
$0\longrightarrow V\longrightarrow M_{n}\longrightarrow \cdots\longrightarrow M_{2}\longrightarrow N_{1}\longrightarrow T_{U}\longrightarrow 0$
By (1), we complete this proof.
\ \hfill $\Box$

\mskip\

By above lemma, we can give the following proposition.

\bg{Pro}\label{pres}
Let $T$ be ${\mathscr X}$-self-orthogonal. Then $\mathrm{Pres}^{n}_{{\mathscr X}}(_{T}{\mathscr X})=\mathrm{Pres}^{n}_{{\mathscr X}}(T)$.

\ed{Pro}

\Pf. we only need to prove that $\mathrm{Pres}^{n}_{{\mathscr X}}(_{T}{\mathscr X})\subseteq\mathrm{Pres}^{n}_{{\mathscr X}}(T)$.
For any $M\in\mathrm{Pres}^{n}_{{\mathscr X}}(_{T}{\mathscr X})$, there is an ${\mathscr X}$-exact sequence
$0\longrightarrow K\longrightarrow N_{n}\longrightarrow \cdots \longrightarrow N_{1}\longrightarrow M\longrightarrow 0$ with
$N_{i}\in$ $_{T}{\mathscr X}$. By Lemma \ref{ltba}, there is an ${\mathscr X}$-exact sequence
$0\longrightarrow V\longrightarrow T_{n}\longrightarrow \cdots \longrightarrow T_{1}\longrightarrow M\longrightarrow 0$ with
$T_{i}\in \mathrm{Add}T$. i.e., $M\in\mathrm{Pres}^{n}_{{\mathscr X}}(T)$.
\ \hfill $\Box$

\mskip\

Next, we give a property of $n$-${\mathscr X}$-tilting module, which is useful to prove the main result of this section (see Theorem \ref{mainresult}).

\bg{Pro}\label{three}
If T is $n$-${\mathscr X}$-tilting, then the following statement are equivalent.

$(1)$ $U\in T^{\mathscr X\bot}$;

$(2)$ $U\in$ $_{T}{\mathscr X}$;

$(3)$ $U\in \mathrm{Pres}^{n}_{\mathscr X}$T.
\ed{Pro}

\Pf. (1)$\Leftrightarrow$ (2) By Lemma \ref{laba} (3).

(2) $\Rightarrow$ (3) Clearly.

(3) $\Rightarrow$ (1) Since $U\in \mathrm{Pres}^{n}_{\mathscr X}T$, there is an $\mathscr X$-exact sequence
$$\xymatrix{0\ar[r]&V\ar[r]&T^{n}\ar[r]&\cdots\ar[r]&T^{1}\ar[r]&U\ar[r]&0}$$
with $T^{i}\in \mathrm{Add}T$ for any $i$. By dimension shifting, we have that
$\mathrm{Ext}_{\mathscr X}^{i+n}(T,~V)\cong \mathrm{Ext}_{\mathscr X}^{i}(T,~U)$. Since ${\mathscr X}$-dim $T \leq n$,
$0=\mathrm{Ext}_{\mathscr X}^{i+n}(T,~V)\cong \mathrm{Ext}_{\mathscr X}^{i}(T,~U)$ for any $i\geq1$ by Proposition 2.4 in \cite{CHEN}.
i.e., $U\in T^{\mathscr X\bot}$.
\ \hfill $\Box$

\mskip\

Recall that a class $\mathcal{C}$ is said to be closed under $n$-$\mathscr X$-images if for any $\mathscr X$-exact sequence
$C_{n}\longrightarrow \cdots \longrightarrow C_{2}\longrightarrow C_{1}\longrightarrow X\longrightarrow 0$ with $C_{i}\in\mathcal{C}$ for all $i$,
then $X\in\mathcal{C}$. It is equivalent to that $\mathrm{Pres}^{n}_{\mathscr X}\mathcal{C}\subseteq \mathcal{C}$.
Obviously, the class $\mathrm{Pres}^{1}_{\mathscr X}T$ is closed under $1$-$\mathscr X$-images. We do
not know that whether or not $\mathrm{Pres}^{n}_{\mathscr X}T$ is closed under $n$-$\mathscr X$-images in general.
But if $\mathrm{Pres}^{n}_{\mathscr X}T$  $=T^{\mathscr X \bot}$, then
$\mathrm{Pres}^{n}_{\mathscr X}T$ is closed under $n$-$\mathscr X$-images, see Lemma \ref{pres=t}.
In particular, if $T$ is $n$-${\mathscr X}$-tilting, $\mathrm{Pres}^{n}_{\mathscr X}T$ is closed under $n$-$\mathscr X$-images by Proposition \ref{three}.

\bg{Lem}\label{pres=t}
Let ${\mathscr A}$ be an Abelian category and $T\in {\mathscr A}$.

$(1)$ If $\mathrm{Pres}^{n}_{\mathscr X}T$  $=T^{\mathscr X \bot}$, then $_{T}{\mathscr X}$  $=T^{\mathscr X \bot}$ is
closed under $n$-$\mathscr X$-images.

$(2)$ ${\mathscr X}$-$dim$ $T \leq n$ if and only if $T^{\mathscr X \bot}$ is
closed under $n$- $\mathscr X$-images.
\ed{Lem}

\Pf. (1) By Lemma \ref{laba} (2), $_{T}{\mathscr X}$  $=T^{\mathscr X \bot}$.
Since $\mathrm{Pres}^{n}_{\mathscr X}(T^{\mathscr X \bot})=$ $\mathrm{Pres}^{n}_{{\mathscr X}}(_{T}{\mathscr X})=\mathrm{Pres}^{n}_{{\mathscr X}}(T)$
(by Proposition \ref{pres}) $=T^{\mathscr X \bot}$, $T^{\mathscr X \bot}$ is closed under $n$-$\mathscr X$-images.

(2) Assume that there is an ${\mathscr X}$-exact sequence
$M_{n}\longrightarrow  \cdots \longrightarrow M_{1}\longrightarrow N\longrightarrow 0$ with
$M_{i}\in T^{\mathscr X \bot}$, set $N_{i}=\mathrm{Ker}(M_{i}\longrightarrow M_{i-1})$ for $1\leq i\leq n$, where $M_{0}=N$.
By dimension shifting, we have that $\mathrm{Ext}^{i+n}_{\mathscr X}(T,~N_{n})\cong\mathrm{Ext}^{i}_{\mathscr X}(T,~N)$.
Since ${\mathscr X}$-$dim$ $T \leq n$, $N\in T^{\mathscr X \bot}$. i.e., $T^{\mathscr X \bot}$ is
closed under $n$- $\mathscr X$-images.

Conversely, for any $L\in {\mathscr A}$. Note that $({\mathscr X}$, ${\mathscr Y})$ is admissible balanced, we have a ${\mathscr Y}$-coresolution of $L$,
$0\longrightarrow L\longrightarrow Y_{0} \longrightarrow \cdots \longrightarrow Y_{n-1}\longrightarrow L_{n}\longrightarrow 0$ with $Y_{i}\in {\mathscr Y}$,
and it is ${\mathscr X}$-exact. Since ${\mathscr Y}\subseteq T^{\mathscr X \bot}$ is closed under $n$- $\mathscr X$-images,
we can obtain that $L_{n}\in T^{\mathscr X \bot}$. By dimension shifting, we have that
$0=\mathrm{Ext}^{i+n}_{\mathscr X}(T,~L)\cong\mathrm{Ext}^{i}_{\mathscr X}(T,~L_{n})$ for any $i$.
So ${\mathscr X}$-$dim$ $T \leq n$ by Lemma 2.4 in \cite{CHEN}.
\ \hfill $\Box$

\mskip\

Finally, we give the main result of this section.

\bg{Th}\label{mainresult}
The following statements are equivalent.

$(1)$ $T$ is $n$-${\mathscr X}$-tilting.

$(2)$ $\mathrm{Pres}^{n}_{\mathscr X}T$  $=T^{\mathscr X \bot}$.
\ed{Th}

\Pf. (1) $\Rightarrow$ (2) By Proposition \ref{three}.

(2) $\Rightarrow$ (1) $(i)$ ${\mathscr X}$-$dim$ $T \leq n$ by Lemma \ref{pres=t}.

$(ii)$ Obviously, $T$ is ${\mathscr X}$-self-orthogonal.

($iii$) For any $X\in\mathscr X$, there is an ${\mathscr X}$-exact sequence
$0\longrightarrow X\longrightarrow Y_{0} \longrightarrow \cdots \longrightarrow Y_{n-1}\longrightarrow Z\longrightarrow 0$ with $Y_{i}\in {\mathscr Y}$
since $({\mathscr X}$, ${\mathscr Y})$ is an admissible balanced pair.
Note that ${\mathscr Y}\subseteq T^{\mathscr X \bot}=$ $_{T}{\mathscr X}$. By Proposition \ref{pres=t}, we have that $Z\in$ $_{T}{\mathscr X}$.
Thus, by Lemma \ref{ltba}, there is an ${\mathscr X}$-exact sequence
$0\longrightarrow A\longrightarrow B\longrightarrow X\longrightarrow 0$ with $A\in$ $_{T}{\mathscr X}$ and $B\in(\widecheck{\mathrm{Add}_{{\mathscr X}}T})_{n}$.
Clearly, it is split as it is ${\mathscr X}$-exact and $X\in\mathscr X$. i.e., $B=A\bigoplus X$.
So, by the dual of the lemma \ref{ltbx}, we have that $X\in(\widecheck{\mathrm{Add}_{{\mathscr X}}T})$. As the definition of the $n$-${\mathscr X}$-tilting module,
$T$ is $n$-${\mathscr X}$-tilting.
\ \hfill $\Box$

\section{Relative tilting modules and Recollement}

Let as Definition \ref{RE} $R(\mathscr A', ~\mathscr A, ~\mathscr A'')$ be a recollement of abelian categories:
$$\xymatrix{
{\mathscr A^\prime} \ar[rrr]|{\ i_{\ast}}&&& {\mathscr A}
\ar@/^/@<2ex>[lll]|{i^{!}}\ar[rrr]|{\
j^{\ast}}\ar@/_/@<-2ex>[lll]|{i^{\ast}} &&& {\mathscr A''}
\ar@/_/@<-2ex>[lll]|{j_{!}}\ar@/^/@<2ex>[lll]|{j_{\ast}} }$$
In this section, we mainly study the relationship of $n$-${\mathscr X}$ tilting modules among three
Abelian categories in a recollement. In a recollement $R(\mathscr A', ~\mathscr A, ~\mathscr A'')$,
we can give a $n$-${\mathscr X}$ tilting module in $\mathscr A$ from two $n$-${\mathscr X}$ tilting modules in $\mathscr A'$ and $\mathscr A''$.
On the contrary, from a $n$-${\mathscr X}$ tilting module in $\mathscr A$, we can also obtain two $n$-${\mathscr X}$ tilting modules
in $\mathscr A'$ and $\mathscr A''$, respectively.
Whenever the abelian category ${\mathscr A}$ (respectively, ${\mathscr A'}$, ${\mathscr A''}$) appears, we always suppose
that it is a part of $R(\mathscr A', ~\mathscr A, ~\mathscr A'')$.

\bg{Th}\label{TH1}%
Suppose that $({\mathscr X}$, ${\mathscr Y})$ is an admissible balanced pair and T is $n$-${\mathscr X}$ tilting module in $\mathscr A$.

$(1)$ If $i^{\ast}$ is an exact functor and $i_{\ast}i^{\ast}(T^{\mathscr X\bot})\subseteq T^{\mathscr X \bot}$,
then $i^{\ast}T$ is $n$-$i^{\ast}{\mathscr X}$ tilting in $\mathscr A'$.

$(2)$ If $j^{\ast}$ is faithful and $j_{\ast}j^{\ast}(T^{\mathscr X\bot})\subseteq T^{\mathscr X \bot}$,
then $j^{\ast}T$ is $n$-$j^{\ast}{\mathscr X}$ tilting in $\mathscr A''$.
\ed{Th}

\Pf. (1) By Theorem \ref{result2}, we know that the pair $(i^{\ast}{\mathscr X}$, $i^{!}{\mathscr Y})=(\mathscr X',~\mathscr Y')$ is admissible balanced.
Note that there is an $\mathscr X$-exact sequence
$0\longrightarrow X_{n}\longrightarrow \cdots\longrightarrow X_{1}\longrightarrow X_{0}\longrightarrow T\longrightarrow0$
with $X_{i}\in \mathscr X$ since $T$ is $n$-${\mathscr X}$ tilting. It is easy to verify that the following exact sequence
$$\xymatrix{
0\ar[r]&i^{\ast}X_{n}\ar[r]&\cdots\ar[r]&i^{\ast}X_{1}\ar[r]&i^{\ast}X_{0}\ar[r]&i^{\ast}T\ar[r]&0
}$$
is a left $\mathscr X'$-resolution of $i^{\ast}T$ as $i^{\ast}$ is exact. i.e., $\mathscr X'$-dim $i^{\ast}T\leq n$.

Since $i^{\ast}$ is exact and $i^{\ast}{\mathscr X}=\mathscr X'$, we have that
$\mathrm{Ext}^{i\geq1}_{\mathscr X'}(i^{\ast}T,~i^{\ast}T^{(I)})\cong \mathrm{Ext}^{i\geq1}_{\mathscr X}(T,~i_{\ast}i^{\ast}T^{(I)})$
for any set $I$. Indeed, let $A^{\bullet}$ be a left $\mathscr X$-resolution of $T$, then $i^{\ast}A^{\bullet}$ is
a left $\mathscr X'$-resolution of $i^{\ast}T$. By definition of the functor $\mathrm{Ext}^{i\geq1}_{\mathscr X}(-,~-)$,
the isomorphism above holds. As $i_{\ast}i^{\ast}(T^{\mathscr X\bot})\subseteq T^{\mathscr X \bot}$, the latter equals zero.
i.e., $i^{\ast}T$ is ${\mathscr X'}$-self-orthogonal.

For any $X\in\mathscr X$, there is an $\mathscr X$-exact sequence
$0\longrightarrow X\longrightarrow T_{0}\longrightarrow \cdots\longrightarrow T_{n}\longrightarrow0$
with $T_{i}\in \mathrm{Add}T$ since $T$ is $n$-${\mathscr X}$ tilting.
Since $i^{\ast}$ is exact, then the following sequence is $\mathscr X'$-exact:
$$\xymatrix{
0\ar[r]&i^{\ast}X\ar[r]&i^{\ast}T_{0}\ar[r]&\cdots\ar[r]&i^{\ast}T_{n}\ar[r]&0
}$$
where $i^{\ast}T_{0}\in \mathrm{Add}i^{\ast}T$.
In conclusion, $i^{\ast}T$ is $n$-$i^{\ast}{\mathscr X}$ tilting.

(2) Similar to (1).
\ \hfill $\Box$

\mskip\

%
%
%
In section 3, we don't know wether the pair $({\mathscr X}$, ${\mathscr Y})$ defined in Theorem \ref{result1} is a balanced
pair or not. But if this is correct, then we have the following conclusion.

\bg{Pro}\label{prop2}%
Suppose that $({\mathscr X'}$, ${\mathscr Y'})$ and $({\mathscr X''}$, ${\mathscr Y''})$ are two admissible balanced pairs in ${\mathscr A'}$ and ${\mathscr A''}$
satisfying the pair $({\mathscr X}$, ${\mathscr Y})$ defined in Theorem \ref{result1} is admissible balanced. If $T'$ (resp., $T''$) is partial $n$-${\mathscr X'}$ tilting
(resp., partial $n$-${\mathscr X''}$ tilting) in $\mathscr A'$ (resp., $\mathscr A''$), and $i^{\ast}$ is exact, then $i_{\ast}T'\bigoplus j_{!}T''$ is partial
$n$-${\mathscr X}$ tilting in $\mathscr A$.
\ed{Pro}

\Pf. (1) As $T'$ is partial $n$-${\mathscr X'}$ tilting, we have the following the ${\mathscr X'}$-exact sequence,
\begin{equation}\label{1}
\begin{split}
\xymatrix{
0\ar[r]&X'_{n}\ar[r]& \cdots\ar[r]&X'_{0}\ar[r]&T'\ar[r]&0
}
\end{split}\tag{4.1}
\end{equation}
with any $X'_{i} \in {\mathscr X'}$. Applying the functor $i_{\ast}$ to (\ref{1}), we have the following sequence,
\begin{equation}\label{2}
\begin{split}
\xymatrix{
0\ar[r]&i_{\ast}X'_{n}\ar[r]& \cdots\ar[r]&i_{\ast}X'_{0}\ar[r]&i_{\ast}T'\ar[r]&0
}
\end{split}\tag{4.2}
\end{equation}
For any $X\in\mathscr X$, applying the functor $\mathrm{Hom}_{\mathscr A}(X,~-)$ to (\ref{2}), we have the following sequence,
\begin{equation}\label{3}
\begin{split}
\xymatrix@C=0.5cm{
0\ar[r]&\mathrm{Hom}_{\mathscr A}(X,~i_{\ast}X'_{n})\ar[r]& \cdots\ar[r]&\mathrm{Hom}_{\mathscr A}(X,~i_{\ast}X'_{0})\ar[r]&\mathrm{Hom}_{\mathscr A}(X,~i_{\ast}T')\ar[r]&0
}
\end{split}\tag{4.3}
\end{equation}
Since $(i^{\ast},~i_{\ast})$ is a adjoint pair, then we have the following commutative diagram.
\begin{equation}\label{4}
\begin{split}
\xymatrix@C=0.5cm{
0\ar[r]&\mathrm{Hom}_{\mathscr A}(X,~i_{\ast}X'_{n})\ar[r]\ar[d]^{\cong}& \cdots\ar[r]&\mathrm{Hom}_{\mathscr A}(X,~i_{\ast}X'_{0})\ar[d]^{\cong}\ar[r]&\mathrm{Hom}_{\mathscr A}(X,~i_{\ast}T')\ar[d]^{\cong}\ar[r]&0\\
0\ar[r]&\mathrm{Hom}_{\mathscr A'}(i^{\ast}X,~X'_{n})\ar[r]& \cdots\ar[r]&\mathrm{Hom}_{\mathscr A'}(i^{\ast}X,~X'_{0})\ar[r]&\mathrm{Hom}_{\mathscr A'}(i^{\ast}X,~T')\ar[r]&0
}
\end{split}\tag{4.4}
\end{equation}
Note that $i^{\ast}{\mathscr X}={\mathscr X'}$ by Theorem \ref{result1}.
From the commutative diagram (\ref{4}), the first row is exact since the sequence (\ref{1}) is ${\mathscr X'}$-exact.
It is easy to verify that $i_{\ast}{\mathscr X'}\subseteq{\mathscr X}$ by the definition of $\mathscr X$.
So the sequence (\ref{2}) is left $\mathscr X$-resolution of $i_{\ast}T'$, i.e., $\mathscr X$- dim $i_{\ast}T'\leq n$.

Note that there is an ${\mathscr X''}$-exact sequence as follows
\begin{equation}\label{5}
\begin{split}
\xymatrix{
0\ar[r]&X''_{n}\ar[r]& \cdots\ar[r]&X''_{0}\ar[r]&T''\ar[r]&0
}
\end{split}\tag{4.5}
\end{equation}
with any $X''_{i} \in {\mathscr X''}$ since partial $T''$ is $n$-${\mathscr X''}$ tilting.
Applying the functor $j_{!}$ to sequence (\ref{5}), we have the following sequence,
\begin{equation}\label{6}
\begin{split}
\xymatrix{
0\ar[r]&j_{!}X''_{n}\ar[r]& \cdots\ar[r]&j_{!}X''_{0}\ar[r]&j_{!}T''\ar[r]&0
}
\end{split}\tag{4.6}
\end{equation}
Further, we can obtain the following sequence,
\begin{equation}\label{7}
\begin{split}
\xymatrix@C=0.5cm{
0\ar[r]&\mathrm{Hom}_{\mathscr A}(j_{!}T'',~Y)\ar[r]&\mathrm{Hom}_{\mathscr A}(j_{!}X''_{0},~Y)\ar[r]& \cdots\ar[r]&\mathrm{Hom}_{\mathscr A}(j_{!}X''_{n},~Y)\ar[r]&0
}
\end{split}\tag{4.7}
\end{equation}
with $Y\in\mathscr Y$. Since $(j_{!},~j^{\ast})$ is a adjoint pair, we can obtain the following commutative diagram.
\begin{equation}\label{8}
\begin{split}
\xymatrix@C=0.5cm{
0\ar[r]&\mathrm{Hom}_{\mathscr A}(j_{!}T'',~Y)\ar[r]\ar[d]^{\cong}&\mathrm{Hom}_{\mathscr A}(j_{!}X''_{0},~Y)\ar[r]\ar[d]^{\cong}& \cdots\ar[r]&\mathrm{Hom}_{\mathscr A}(j_{!}X''_{n},~Y)\ar[r]\ar[d]^{\cong}&0\\
0\ar[r]&\mathrm{Hom}_{\mathscr A''}(T'',~j^{\ast}Y)\ar[r]&\mathrm{Hom}_{\mathscr A''}(X''_{0},~j^{\ast}Y)\ar[r]& \cdots\ar[r]&\mathrm{Hom}_{\mathscr A''}(X''_{n},~j^{\ast}Y)\ar[r]&0
}
\end{split}\tag{4.8}
\end{equation}
It is easily see that $j_{!}\mathscr X''\subseteq \mathscr X$ by the definition of $\mathscr X$ and Proposition \ref{propR} (4).
Note that $j^{\ast}{\mathscr Y}={\mathscr Y''}$ by Theorem \ref{result1}. Thus the second row in the commutative diagram (\ref{8}) is exact since
$({\mathscr X''}$, ${\mathscr Y''})$ is an admissible balanced pair.
So the sequence (\ref{6}) is left $\mathscr X$-resolution of $j_{!}T''$, i.e., $\mathscr X$- dim $j_{!}T''\leq n$.
Consequently, $\mathscr X$- dim $j_{!}T''\bigoplus i_{\ast}T'\leq n$.

(2) $\mathrm{Ext}^{i\geq1}_{\mathscr X}(i_{\ast}T',~(i_{\ast}T')^{(I)})\cong \mathrm{Ext}^{i\geq1}_{\mathscr X'}(i^{\ast}i_{\ast}T',~T'^{(I)})
\cong \mathrm{Ext}^{i\geq1}_{\mathscr X'}(T',~T'^{(I)})=0$. The first isomorphism holds since $(i^{\ast},~i_{\ast})$ is an adjoint pair,
they are exact functor and $i^{\ast}\mathscr X=\mathscr X'$. The second isomorphism holds as $i^{\ast}i_{\ast}=1_{\mathscr A'}$ by
Proposition \ref{propR}.

$\mathrm{Ext}^{i\geq1}_{\mathscr X}(j_{!}T'',~(j_{!}T'')^{(I)})\cong \mathrm{Ext}^{i\geq1}_{\mathscr X''}(T'',~j^{\ast}j_{!}T''^{(I)})
\cong \mathrm{Ext}^{i\geq1}_{\mathscr X''}(T'',~T''^{(I)})=0$. Note that $j_{!}$ is also exact by Proposition \ref{propR} (4).
Thus the first isomorphism holds.

$\mathrm{Ext}^{i\geq1}_{\mathscr X}(i_{\ast}T',~(j_{!}T'')^{(I)})\cong \mathrm{Ext}^{i\geq1}_{\mathscr X'}(T',~i^{!}j_{!}T''^{(I)})=0$ by Proposition \ref{propR} (4).

$\mathrm{Ext}^{i\geq1}_{\mathscr X}(j_{!}T'',~(i_{\ast}T')^{(I)})\cong \mathrm{Ext}^{i\geq1}_{\mathscr X''}(T'',~j^{\ast}i_{\ast}T'^{(I)})=0$.
Thus $i_{\ast}T'\bigoplus j_{!}T''$ is ${\mathscr X}$-self-orthogonal.
%

In conclusion, $i_{\ast}T'\bigoplus j_{!}T''$ is partial $n$-${\mathscr X}$ tilting in $\mathscr A$.
\ \hfill $\Box$

\mskip\

In fact, from the above proof, we have that both $i_{\ast}T'$ and $j_{!}T''$ are partial $n$-${\mathscr X}$ tilting in $\mathscr A$.

\bg{Th}\label{TH2}
Under these conditions in Theorem \ref{result1}, if $T''$ is $n$-${\mathscr X''}$ tilting in $\mathscr A''$ and $j_{!}$ is exact, then $j_{!}T''$ is $n$-${\mathscr U}$ tilting in $\mathscr A$,
where $({\mathscr U}$, ${\mathscr V})$ = $(j_{!}j^{\ast}{\mathscr X}$, $j_{\ast}j^{\ast}{\mathscr Y})$.
\ed{Th}

\Pf. It follows from Theorem \ref{result1} that $(j_{!}j^{\ast}{\mathscr X}$, $j_{\ast}j^{\ast}{\mathscr Y})$ is admissible balanced.

(1) Since $T''$ is $n$-${\mathscr X''}$ tilting there is an ${\mathscr X''}$-exact sequence as follows
\begin{equation}\label{9}
\begin{split}
\xymatrix{
0\ar[r]&X''_{n}\ar[r]& \cdots\ar[r]&X''_{0}\ar[r]&T''\ar[r]&0
}
\end{split}\tag{4.9}
\end{equation}
with any $X''_{i} \in {\mathscr X''}$. Applying the functor $j_{!}$ to sequence (\ref{9}), we have the following sequence,
\begin{equation}\label{10}
\begin{split}
\xymatrix{
0\ar[r]&j_{!}X''_{n}\ar[r]& \cdots\ar[r]&j_{!}X''_{0}\ar[r]&j_{!}T''\ar[r]&0
}
\end{split}\tag{4.10}
\end{equation}
For any $X\in\mathscr X$, we can obtain the following commutative diagram since $j_{!}$ is fully faithful.
\begin{equation}\label{11}
\begin{split}
\xymatrix@C=0.5cm{
0\ar[r]&\mathrm{Hom}_{\mathscr A}(j_{!}j^{\ast}X,~j_{!}X''_{n})\ar[r]\ar[d]^{\cong}& \cdots\ar[r]&\mathrm{Hom}_{\mathscr A}(j_{!}j^{\ast}X,~j_{!}X''_{0})\ar[d]^{\cong}\ar[r]&\mathrm{Hom}_{\mathscr A}(j_{!}j^{\ast}X,~j_{!}T'')\ar[d]^{\cong}\ar[r]&0\\
0\ar[r]&\mathrm{Hom}_{\mathscr A'}(j^{\ast}X,~X''_{n})\ar[r]& \cdots\ar[r]&\mathrm{Hom}_{\mathscr A'}(j^{\ast}X,~X''_{0})\ar[r]&\mathrm{Hom}_{\mathscr A'}(j^{\ast}X,~T'')\ar[r]&0
}
\end{split}\tag{4.11}
\end{equation}
Note that $j^{\ast}{\mathscr X}={\mathscr X''}$. From the sequence (\ref{9}), we know that the second row in (\ref{11}) is exact, thus , the first row is also exact.
Consequently, $\mathscr U$- dim $j_{!}T''\leq n$ from the sequence (\ref{10}).

(2) Similar to the proof of Proposition \ref{prop2}, $j_{!}T''$ is ${\mathscr U}$-self-orthogonal.

(3) For any $X\in\mathscr X$, $j^{\ast}X\in\mathscr X''$. Then there is a ${\mathscr X''}$-exact sequence
$$\xymatrix{
0\ar[r]&j^{\ast}X\ar[r]&T''_{0}\ar[r]&\cdots\ar[r]&T''_{n}\ar[r]&0
}
$$
with $T''_{i}\in \mathrm{Add}T''$ since $T''$ is $n$-${\mathscr X''}$ tilting. Applying the functor $j_{!}$ to the above sequence, we have the following sequence,
\begin{equation}\label{12}
\begin{split}
\xymatrix{
0\ar[r]&j_{!}j^{\ast}X\ar[r]&j_{!}T''_{0}\ar[r]&\cdots\ar[r]&j_{!}T''_{n}\ar[r]&0
}
\end{split}\tag{4.12}
\end{equation}
Similar to the discuss of (1), the sequence (\ref{12}) is our desired.

To sum up, we claim that $j_{!}T''$ is $n$-${\mathscr U}$ tilting.

{\small

}

\end{document}